# A Radial Basis Function Partition of Unity Method for Steady Flow Simulations

Francisco Bernal,[*] Ali Safdari-Vaighani,[†] Elisabeth Larsson[‡]

**Abstract**

A methodology is presented for the numerical solution of nonlinear elliptic systems in unbounded domains, consisting of three elements. First, the problem is posed on a finite domain by means of a proper nonlinear change of variables. The compressed domain is then discretised, regardless of its final shape, via the radial basis function partition of unity method. Finally, the system of nonlinear algebraic collocation equations is solved with the trust-region algorithm, taking advantage of analytically derived Jacobians. We validate the methodology on a benchmark of computational fluid mechanics: the steady viscous flow past a circular cylinder. The resulting flow characteristics compare very well with the literature. Then, we stress-test the methodology on less smooth obstacles—rounded and sharp square cylinders. As expected, in the latter scenario the solution is polluted by spurious oscillations, owing to the presence of boundary singularities.

## 1 Introduction

We consider the numerical solution of systems of nonlinear elliptic boundary value problems (BVPs) in unbounded domains. Examples of applications modeled by such BVPs are the electrostatic potential of a biomolecule in a solvent [32] (by the Poisson-Boltzmann equation); the exterior wave scattering problem with state-dependent permittivity [9] (by the nonlinear Helmholtz equation); or the viscous flow past immersed bodies, which models fluid-structure interaction. This paper is devoted to the latter.

When solving those problems, boundary conditions (BCs) are typically enforced on the boundary of a finite computational domain. Since such a truncation pollutes the solution in the region of interest, the computational boundary must be set far out from it. Upon discretisation, this results in large and taxing nonlinear systems of equations.

[*]Department of Mathematics, Universidad Carlos III de Madrid, Spain (fcoberna@math.uc3m.es)
[†]Department of Mathematics, Allameh Tabataba'i University, Iran (asafdari@atu.ac.ir)
[‡]Department of Information Technology, Uppsala University, Sweden (elisabeth.larsson@it.uu.se)

In this paper (restricted to two-dimensional geometries for the sake of simplicity) we have followed the far less common alternative of compressing $\mathbb{R}^2$ into a bounded domain by a proper transformation—thus enabling the actual BCs at infinity to be enforced—see Figure 1. In general, not only will such a transformation lead to a scattered discretisation of the transformed domain, but it will spoil whatever weak (or otherwise convenient form) the equations may have had, in the first place. For those reasons, we have chosen to solve the transformed BVP by radial basis function (RBF) collocation [10], which is ideally suited to solving partial differential equations (PDEs) in the strong sense on unstructured nodes inside irregular domains.

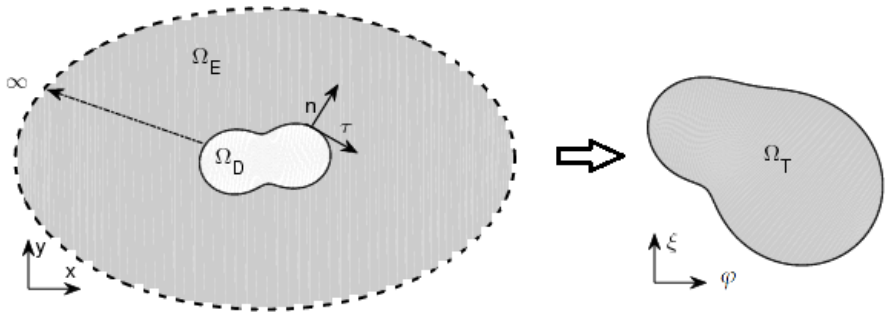

Figure 1: The body $\Omega_D$ and the surrounding unbounded physical domain $\Omega_E$ (where the PDEs are to be solved), mapped onto the bounded computational domain $\Omega_T$. (The unit normal and tangent vector are depicted on one point of $\partial \Omega_D$.)

Specifically, we have tested our methodology on the steady flow of a viscous fluid past a circular cylinder—a benchmark for numerical methods in computational fluid mechanics [5, 11, 12, 13, 22, 24, 27, 34]. (We shall describe it in a moment.) In order to solve the (transformed) steady-state Navier-Stokes equations on the transformed, bounded domain, we have used the RBF - partition of unity (RBF-PU) method [30]. RBF-PU combines Kansa's RBF collocation on overlapping patches with a partition of unity which glues them together, as described in [1]. (A preliminary version without results was given in [4].) It retains spectral accuracy while giving rise to a sparse—rather than full—collocation matrix, which is much better conditioned, and thus capable of accommodating significantly more collocation equations (see also [21]).

The physical setup consists in a long, circular cylinder (along the $z$-axis) immersed in a viscous fluid which flows steadily along the $x$-direction; the problem can thus be regarded as two dimensional—see Figure 2 (left). Let $U_\infty$ be the velocity of the fluid at infinity, $\nu$ its kinematic viscosity, $2R$ the cylinder diameter, and let $Re = 2RU_\infty/\nu$ be the (diameter-based) Reynolds number. Far away from the cylinder, the flow is unaffected by it; whereas immediately past the cylinder, an elongated area of flow recirculation occurs after $Re \gtrsim 5$, called the wake. For low $Re$, the flow is steady and symmetric along the $x$-axis, and can be solved by substituting a zero for the time derivative in the (two-dimensional) Navier-Stokes equations, yielding—in the velocity-pressure, or natural variables, formulation—a nonlinear system of three equations. However, the actual flow becomes unstable after $Re \gtrsim 47$ [16], when the phenomenon of vortex shedding kicks in, whereby eddies detach alternatingly from either side of the cylinder surface. Past this threshold the flow is oscillating, rather than steady, and can no longer be modeled by time-independent equations. Consequently, all of our numerical experiments took place below the instability threshold.

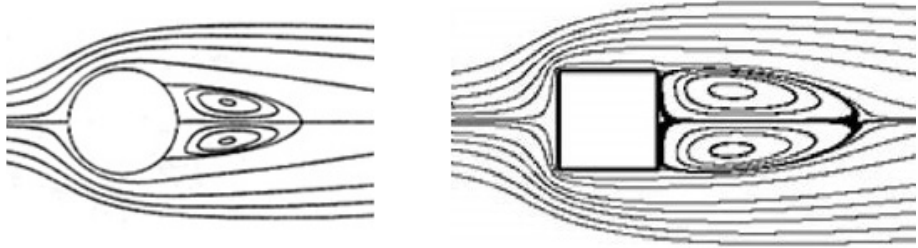

Figure 2: Steady viscous flow past a circular (left) and square (right) cylinder. Note the recirculation region, or bubble, and the velocity point singularities.

It was shown in [3] that RBF discretisations of nonlinear elliptic BVPs are best tackled with the trust-region algorithm, owing to the accuracy of the Jacobians—calculated via the Fréchet derivative of the nonlinear operator—and to the simplicity in constructing and using them. RBF-PU Jacobians are, in addition to this, sparse. In all of the experiments carried out in this work, the basic dogleg trust-region subproblem approximation proved adequate to find the solution within a few iterations, even near the threshold $Re$ value.

To the best of our knowledge, this work actually is the first time that RBF-PU has been applied to a nonlinear elliptic PDE. In fact, examples of nonlinear elliptic PDEs solved with RBFs are conspicuously scarce in the literature; let us briefly review some of them (besides [3]). Paper [14] is concerned with semilinear Poisson equations, based on stencil-based Hermite approximations and Picard iterations for the right hand side, through linearisation. Nonlinear mechanics—discretised by Kansa's method with and without Hermite—are addressed in [35] by means of incremental iteration. RBF-FD (finite differences) are employed in [15] with a nonlinear material, thanks to a direct iterative solution process. The RBF-FD PhD thesis [26] on nonlinear mechanics is perhaps the closest to our approach. It uses Newton's method with exact Jacobian for thermoelasticity. For a related elasto-viscoplastic solid mechanics problem, the Jacobian is exact for the linear part and numerical for the non-linear part. In all of the previous examples, the PDEs were handled in strong form. In contrast, laminar flow is resolved in [17] by letting Matlab's `fsolve` and `lsqnonlin` find the Kansa's weights in oversampled form (and the shape parameter as well). As usual, the Jacobian/gradient is not provided to the solver, though.

Despite the fact that our transformed RBF-PU discretisation stretches from the body surface to infinity, it adequately reproduces the flow patterns in the wake. Our results for the drag coefficient and other standard quality metrics, such as the bubble shape, are in good agreement with the literature. They also attest to residual and error convergence as the discretisation is refined. In sum, our method yields competitive re-sults for the benchmark circular cylinder, while avoiding the ambiguity in (and the cost entailed by) the location of the BCs on a distant computational boundary. Since our formulation solves for the natural variables, extending it to three dimensions would be straightforward (by adding an equation for the vertical velocity to the system). Moreover, it is meshless and amenable to the strong form of the PDE and BCs.

On the other hand, spectral methods typically struggle with insufficiently smooth problems. In the second battery of experiments, we have stress-tested our methodology by considering edgier obstacles to the flow than a circular cylinder. Concretely, we replace the latter by a rounded square cross-section, and gradually sharpen it up into a perfectly square cylinder—see Figure 2 (right). There, the flow velocity undergoes abrupt changes on several isolated points of the immersed bodies—emphatically so, on the upstream corners of the square cylinder.

As expected, our numerical results deteriorate during that transition—even though most quality metrics do roughly keep up with the references. While we have not managed to satisfactorily fix this challenging issue, we do believe that those results are nonetheless informative and worth reporting, since nicely smooth problems are more common in the RBF literature than they are in realistic applications.

Following this introduction, the remainder of the paper is organised as follows. In Section 2, we introduce the nonlinear transformation of the half plane outside of the circular cylinder into a bounded domain, and the resulting transformed time-independent Navier-Stokes equations for velocity and pressure. In Section 3, we explain how to solve them with the RBF-PU method and the trust-region algorithm. In Section 4, the steady flow past a circular cylinder is numerically studied with said methodology. In order to test the limits of the latter, in Section 5 the flows past rounded and square cylinders are tackled as well. Section 6 is a concluding commentary on the strengths and limitations of the method, and an invitation for future work.

## 2 Transformed Navier-Stokes equations

A fixed, infinite circular cylinder of radius $R$ is immersed in a fluid of kinematic viscosity $\nu > 0$, which flows steadily perpendicularly to the cylinder with far-field velocity $U_\infty > 0$ and far-field pressure $P_\infty$. The Reynolds number (based on the diameter) is $Re = 2RU_\infty/\nu$. By symmetry, the problem is two-dimensional, with the obstacle being the circular cross section. Let $(x,y)$ be a Cartesian dimensionless frame (with $R = 1$) centred at the axis; $(r, \varphi)$ polar coordinates (where $\varphi = 0$ marks the direction of advance of the flow); $\mathbf{u} = (u, v)$ the dimensionless velocity field (with $U_\infty = 1$); $P$ the pressure and $p = P - P_\infty$. The steady Navier-Stokes equations are then given by

$$\frac{Re}{2}(\mathbf{u}\cdot\nabla)\mathbf{u} = \nabla^2\mathbf{u} - \nabla p, \qquad \nabla\cdot\mathbf{u} = 0. \tag{1}$$

They are supplemented with five boundary conditions:

$$\begin{aligned}
u(r=\infty,\varphi) &= 1, & v(r=\infty,\varphi) &= 0, & &\text{(far-field flow)} \\
u(r=1,\varphi) &= 0, & v(r=1,\varphi) &= 0, & &\text{(non-slip condition)} \\
& & p(r=\infty,\varphi) &= 0. & &\text{(far-field pressure)}
\end{aligned} \tag{2}$$

Rather than taking a large finite domain and enforcing the BCs far away from the cylinder, the infinite domain is compressed into the rectangle $[0,\ell] \times [0,\pi]$ via the following transformation[1]:

$$\xi = \ell(1 - 1/r) \qquad (\text{ such that } \xi(r=1) = 0 \text{ and } \xi(r=\infty) = \ell\ ). \tag{3}$$

(The stretching factor $\ell \geq 1$ will later be handy for discretisation purposes.) Moreover, $\partial/\partial r = (1/\ell)(\ell - \xi)^2 \partial/\partial \xi$, so that the unit vectors point in the same direction: $\mathbf{i}_r = \mathbf{i}_\xi$. Denoting as $\dot\xi$ and $\dot\varphi$ the components of the fluid velocity in the new coordinates, the velocity field is transformed as

$$\mathbf{u} = \dot\xi\,\mathbf{i}_r + (1 - \xi/\ell)\dot\varphi\,\mathbf{i}_\varphi. \tag{4}$$

(The dot notation is intuitive for velocities and useful so as not to overload the notation, but recall that the problem is steady and involves no time derivatives.) Let us use the notation $\partial^2_{\xi\varphi}\dot\varphi = \partial^2\dot\varphi/\partial\xi\partial\varphi$, $\partial_\varphi p = \partial p/\partial\varphi$, etc. The steady, two dimensional Navier-Stokes equations in the variables $(\xi,\varphi)$ are transformed into

$$\mathscr{W}_1 = \mathscr{W}_2 = \mathscr{W}_3 = 0, \tag{5}$$

where:

$$\mathscr{W}_1 = \frac{Re}{2}\left[(\ell-\xi)\dot\xi\partial_\xi\dot\xi + \dot\varphi\partial_\varphi\dot\xi - \dot\varphi^2 + (\ell-\xi)\partial_\xi p\right] \tag{6}$$

$$- \frac{1}{\ell}(\ell-\xi)^3 \partial^2_{\xi\xi}\dot\xi - \frac{1}{\ell}(\ell-\xi)\partial^2_{\varphi\varphi}\dot\xi + \frac{1}{\ell}(\ell-\xi)^2\partial_\xi\dot\xi + \frac{2}{\ell}(\ell-\xi)\partial_\varphi\dot\varphi + \frac{1}{\ell}(\ell-\xi)\dot\xi,$$

$$\mathscr{W}_2 = \frac{Re}{2}\left[(\ell-\xi)\dot\xi\partial_\xi\dot\varphi + \dot\varphi\partial_\varphi\dot\varphi + \dot\xi\dot\varphi + \partial_\varphi p\right] \tag{7}$$

$$- \frac{1}{\ell}(\ell-\xi)^3\partial^2_{\xi\xi}\dot\varphi - \frac{1}{\ell}(\ell-\xi)\partial^2_{\varphi\varphi}\dot\varphi + \frac{1}{\ell}(\ell-\xi)^2\partial_\xi\dot\varphi - \frac{2}{\ell}(\ell-\xi)\partial_\varphi\dot\xi + \frac{1}{\ell}(\ell-\xi)\dot\varphi,$$

$$\mathscr{W}_3 = (\ell-\xi)\partial_\xi\dot\xi + \partial_\varphi\dot\varphi + \dot\xi. \tag{8}$$

Moreover, since the problem is symmetric along the $x$ axis, only $0 \leq \varphi \leq \pi$ has to be considered. (The boundary conditions along the $x$ axis are now of reflecting type to enforce the symmetry.) The transformed BCs are sketched in Figure 3.

**Extension to other geometries.** Since the Navier-Stokes equations are posed for velocity-pressure, it is straightforward to extend them to three dimensions, by adding one fourth nonlinear equation $\mathscr{W}_4$. (Unlike in the streamfunction-vorticity formulation, which has one less equation, but is restricted to planar flows.)

While the presented transformation has been tailored to the unbounded domain outside a (half) circle, it can be adapted to other obstacle cross-sections $\Gamma$ as long as its surface can be described as $\Gamma = \partial\Omega_T : f(\xi,\varphi) = 0$. (For instance, if $\Omega_T$ is convex, then $\Gamma$ can be described in polar coordinates.) This is the case of the steady flows past a square obstacle with rounded or sharp corners, discussed in Section 5. In those cases, the transformed domain is no longer a rectangle in the $(\xi,\varphi)$ variables. As long as the transformed boundary is smooth, this will not be an issue, thanks to the geometric flexibility of the RBF-PU method, which is discussed next.

---

[1]This mapping appears in [25], but beware of a mistake or typo in the last two terms of formula (7) there.

$$\dot{\xi}(\xi=\ell,\varphi)=\cos\varphi\,,\ \dot{\varphi}(\xi=\ell,\varphi)=-\sin\varphi\,,\ p(\xi=\ell,\varphi)=0$$

$$\dot{\xi}(\xi=0,\varphi)=0$$
$$\dot{\varphi}(\xi=0,\varphi)=0$$

$$\frac{\partial\dot{\xi}}{\partial\varphi}(\xi,\varphi=\pi)=0 \qquad \frac{\partial\dot{\xi}}{\partial\varphi}(\xi,\varphi=0)=0$$
$$\frac{\partial\dot{\varphi}}{\partial\varphi}(\xi,\varphi=\pi)=0 \qquad \frac{\partial\dot{\varphi}}{\partial\varphi}(\xi,\varphi=0)=0$$
$$\frac{\partial p}{\partial\varphi}(\xi,\varphi=\pi)=0 \qquad \frac{\partial p}{\partial\varphi}(\xi,\varphi=0)=0$$

Figure 3: Transformed BCs for the circular cylinder. The infinite annular section ($1 \leq r < \infty, 0 \leq \varphi \leq \pi$) sketched above has been compressed into the rectangle $[0,\ell] \times [0,\pi]$.

## 3 Methodology

### 3.1 Overview of Kansa's method

Let $\mathbf{q} = (\xi, \varphi) \in \Omega_T \subset \mathbb{R}^2$, and let the *pointset* $\{\mathbf{q}_1, \ldots, \mathbf{q}_N\}$ be a discretisation of $\Omega_T$ and its boundary $\partial \Omega_T$ into $N$ scattered, distinct points (called *nodes*). The solution components $\{\dot{\xi}(\mathbf{q}), \dot{\varphi}(\mathbf{q}), p(\mathbf{q})\}$ is approximated by *RBF interpolants*. For instance,

$$p(\xi, \varphi) = p(\mathbf{q}) = \sum_{i=1}^{N} \lambda_i \phi_i(\|\mathbf{q}\|). \tag{9}$$

*Remark.* For notational convenience, we use the symbols $\dot{\xi}, \dot{\varphi}$, and $p$ both for the exact solution of equations (5) and for their RBF interpolants.

Above, $\|\cdot\|$ is the Euclidean norm, and $\phi_i(\mathbf{q})$ is the chosen RBF, which also contains a *shape parameter* $\varepsilon > 0$. For example, the inverse multiquadric RBF is

$$\phi_i(\mathbf{q}) = \frac{1}{\sqrt{1 + (\varepsilon \|\mathbf{q} - \mathbf{q}_i\|)^2}}. \tag{10}$$

The RBF coefficients $\lambda_1, \ldots, \lambda_N$ can be found by collocation. Enforcing that $p(\mathbf{q})$ in (9) interpolates the nodal pressures $p(\mathbf{q}_1), \ldots, p(\mathbf{q}_N)$ leads to

$$\begin{pmatrix} p_1 \\ \vdots \\ p_N \end{pmatrix} = \begin{pmatrix} \phi_1(\mathbf{q}_1) & \cdots & \phi_N(\mathbf{q}_1) \\ \vdots & \ddots & \vdots \\ \phi_1(\mathbf{q}_N) & \cdots & \phi_N(\mathbf{q}_N) \end{pmatrix} \begin{pmatrix} \lambda_1 \\ \vdots \\ \lambda_N \end{pmatrix} \Rightarrow \vec{\lambda} = [\phi]^{-1} \vec{p}, \tag{11}$$

where we have introduced the notation of *nodal vectors and matrices*, which we extend to functions: given $f : \Omega \mapsto \mathbb{R}$ and a set of collocation points $Q = \{\mathbf{q}_1, \ldots, \mathbf{q}_{\#Q}\}$, then $\vec{f}(Q) \in \mathbb{R}^{\#Q} = [f(\mathbf{q}_1), \ldots, f(\mathbf{q}_{\#Q})]$. Usually, we will just write $\vec{f}$ instead of $\vec{f}(Q)$, so that the size of the nodal vector must be inferred from the context.

Note that (with a fixed $\varepsilon$) $[\phi]$ above is symmetric and, if $\phi$ is a positive definite RBF, invertible. This allows one to express the RBF interpolant in terms of the (unknown) nodal pressures rather than RBF coefficients, by taking each row of $[\phi]^{-1}$, as

$$p(\mathbf{q}) = \vec{\phi}^T(\mathbf{q})\vec{\lambda} = \vec{\phi}^T(\mathbf{q})[\phi]^{-1}\vec{p} = \vec{\psi}^T(\mathbf{q})\vec{p}, \tag{12}$$

where

$$\vec{\psi}(\mathbf{q}) = [\phi]^{-1}\vec{\phi}(\mathbf{q}) \Rightarrow \psi_i(\mathbf{q}_j) = \delta_{ij} \text{ (Kronecker's delta)}. \tag{13}$$

The functions $\psi_i(\mathbf{q})$ are the *cardinal basis functions* of the RBF $\phi$ and the pointset.

Linear BVPs can readily be solved as follows. Let the PDE be defined by the interior operator $\mathscr{L}^{PDE} p = f$ and the BCs by the boundary operator $\mathscr{L}^{BC} p = g$. For notational convenience, let the BVP be described by $\mathscr{L} = h$, where $h(\mathbf{q})$ and $\mathscr{L}$ are

$$\mathscr{L} = \begin{cases} \mathscr{L}^{PDE}, & \text{if } \mathbf{q} \in \Omega_T, \\ \mathscr{L}^{BC}, & \text{if } \mathbf{q} \in \partial\Omega_T, \end{cases} \quad h(\mathbf{q}) = \begin{cases} f(\mathbf{q}), & \text{if } \mathbf{q} \in \Omega_T, \\ g(\mathbf{q}), & \text{if } \mathbf{q} \in \partial\Omega_T. \end{cases} \tag{14}$$

Applying $\mathscr{L}$ on the RBF interpolant of $p$ yields the square linear system

$$\mathscr{L} p(\mathbf{q}) = \left[\mathscr{L}\phi_1(\mathbf{q}), \ldots, \mathscr{L}\phi_N(\mathbf{q})\right]^T [\phi]^{-1} \vec{p} = \vec{h}, \tag{15}$$

whose solution is $\vec{p}$, and $p(\mathbf{q})$ can be reconstructed throughout $\Omega_T \cup \partial\Omega_T$ by (12).

The just described algorithm is easy to code, meshfree, geometrically flexible and spectrally convergent for smooth problems. A drawback is that the last property comes at the expense of a fully populated matrix in (15). For this reason, it is often thought that Kansa's method loses many of its advantages after a couple thousand nodes. The RBF-PU method pushes $N$ significantly beyond that without loss of performance.

### 3.2 Radial Basis Function - Partition of Unity (RBF-PU) method

Let $\Omega_T \subset \mathbb{R}^2$ be a bounded domain and $\{\Omega_i\}_{i=1}^P$ be an open *cover* of $\Omega$ with some mildly overlapping *patches*: $\Omega_T \subset \{\Omega_i\}_{i=1}^P$. The amount of overlap between patches must be limited such that at most $K$ patches overlap at any given point, i.e.

$$\Xi(\mathbf{q}) = \{k | \mathbf{q} \in \Omega_k\}, \quad |\Xi(\mathbf{q})| \leq K, \quad \forall \mathbf{q} \in \Omega_T. \tag{16}$$

Throughout this paper, we choose to define the patches as discs and ellipses, whereby we can guarantee that the domain is covered, and regulate the amount of overlap. Associated with these patches, we construct a *partition of unity*, i.e. a set of nonnegative, compactly supported, continuous *weight functions* $w_i : \Omega_i \to \mathbb{R}$, $i = 1, \ldots, P$ such that

$$\sum_{i=1}^P w_i(\mathbf{q}) = 1. \tag{17}$$

Weight functions can be constructed by Shepard's method [33] as follows:

$$w_i(\mathbf{q}) = \frac{\Phi_i(\mathbf{q})}{\sum_{k=1}^P \Phi_k(\mathbf{q})}, \quad i = 1, \ldots, P, \tag{18}$$

where $\Phi_i(\mathbf{q}) = \Phi\left(\frac{\|\mathbf{q}-\mathbf{c}_i\|}{\rho_i}\right)$ are compactly supported radial functions on $\Omega_i$ with centre $\mathbf{c}_i$ and radius $\rho_i$. $\Phi$ must be smooth enough to support the differential operators of the problem to be solved. For instance, the $C^2$ Wendland function [36] is

$$\Phi^W(r) = (4r+1)(1-r^4)^+. \qquad (19)$$

It is sometimes convenient to use elliptic patches [30]. For instance, an elliptic patch aligned with the axes centred at $(c_{i,1}, c_{i,2})$ with semiaxes $(\rho_{i,1}, \rho_{i,2})$ is given by

$$\Phi_i(q_1, q_2) = \Phi^W\left(\sqrt{\frac{(q_1-c_{i,1})^2}{\rho_{i,1}^2} + \frac{(q_2-c_{i,2})^2}{\rho_{i,2}^2}}\right). \qquad (20)$$

Let us take the pressure as illustration. A global partition of unity approximation to $p(\mathbf{q})$ over $\overline{\Omega}_T$ is formed as a weighted sum of local RBF approximations $p_j(\mathbf{q})$ on each patch $\Omega_j$ of the cover, "glued together" by means of the partition of unity:

$$p(\mathbf{q}) = \sum_{j \in \Xi(\mathbf{q})} w_j(\mathbf{q}) p_j(\mathbf{q}). \qquad (21)$$

Expressing $p_j(\mathbf{q})$ in terms of the nodal values on patch $j$, and by linearity,

$$\mathscr{L} p(\mathbf{q}) = \sum_{j \in \Xi(\mathbf{q})} \sum_{k \in \Omega_j} \mathscr{L}\left(w_j(\mathbf{q}) \psi_k(\mathbf{q})\right) p_k. \qquad (22)$$

Specifically, partial derivatives can be computed according to the formula

$$\frac{\partial^{|\alpha|}}{\partial q^\alpha} p(\mathbf{q}) = \sum_{j \in \Xi(\mathbf{q})} \sum_{k \in \Omega_j} \frac{\partial^{|\alpha|}}{\partial q^\alpha}\left[w_j(\mathbf{q}) \psi_k(\mathbf{q})\right] p_k = \sum_{j \in \Xi(\mathbf{q})} \sum_{k \in \Omega_j} \left[\sum_{\beta \leq \alpha} \binom{\alpha}{\beta} \frac{\partial^{|\alpha-\beta|} w_j}{\partial q^{\alpha-\beta}} \frac{\partial^{|\beta|} \psi_k}{\partial q^\beta}\right] p_k, \qquad (23)$$

where $\partial^{|\alpha|}/\partial^\alpha$ is the usual multi-index notation [30]. For instance, the angular derivative at a node with $\xi = \xi'$ and $\varphi = \varphi'$ is

$$\partial_\varphi p(\xi', \varphi') = \sum_{j \in \Xi(\mathbf{q}')} \sum_{k \in \Omega_j} \left[\frac{\partial w_j}{\partial q_y}(\mathbf{q}') \psi_k(\mathbf{q}') + w_j(\mathbf{q}') \frac{\partial \psi_k}{\partial q_y}(\mathbf{q}')\right] p_k. \qquad (24)$$

We stress the fact that partial derivatives of an RBF interpolant is a linear combination of its nodal values, with the coefficients depending only on the discretisation (i.e. pointset, cover, partition of unity and RBFs), but independent of the function being differentiated. Therefore, one can explicitly compute the matrices for the nodal vectors of the required derivatives at start, and reuse them as needed. (For instance, in a nonlinear loop.) In the previous example, calling $\left((\partial_\varphi p)_1, \ldots, (\partial_\varphi p)_N\right)^T =: \partial_\varphi \vec{p}$,

$$\partial_\varphi \vec{p} =: [\partial_\varphi] \vec{p}, \text{ with } [\partial_\varphi]_{mn} = \begin{cases} \sum_{j \in \Xi(\mathbf{q}_n)} \left[\frac{\partial w_j}{\partial \varphi}(\mathbf{q}_m) \psi_n(\mathbf{q}_m) + w_j(\mathbf{q}_m) \frac{\partial \psi_n}{\partial \varphi}(\mathbf{q}_m)\right] & \text{if } \Xi(\mathbf{q}_m) \cap \Xi(\mathbf{q}_m) \neq \emptyset, \\ 0 & \text{otherwise.} \end{cases} \qquad (25)$$

Note that matrices such as $[\partial_\varphi]$ are sparse because only entries with indices associated to nodes in overlapping patches are nonzero.

Later in this section, we shall need to evaluate derivatives on restricted sets of collocation nodes $Q \subset \{1,\ldots,N\}$ with #$Q$ ordered nodes. Let $1 \leq Q(1) < Q(2) < \ldots < Q(\#Q) \leq N$ be the indices of the nodes in $Q$, and, in the previous example, let $[\partial_\varphi]_Q$ be the rectangular matrix of size $\#Q \times N$ defined as

$$[\partial_\varphi]_Q := \text{submatrix formed by the rows } Q(1),\ldots,Q(\#Q) \text{ of } [\partial\varphi]. \tag{26}$$

For the identity operator in particular,

$$([1]_Q)_{ij} = \begin{cases} 1 \text{ if } Q(i) = j, \\ 0 \text{ otherwise.} \end{cases} \tag{27}$$

Analogously, let us define the RBF interpolants of $\dot{\xi}(\xi,\varphi)$ and $\dot{\varphi}(\xi,\varphi)$ as

$$\dot{\xi}(\mathbf{q}) = \sum_{j \in \Xi(\mathbf{q})} \sum_{k \in \Omega_j} \left(w_j(\mathbf{q})\psi_k(\mathbf{q})\right) \dot{\xi}_k, \qquad \dot{\varphi}(\mathbf{q}) = \sum_{j \in \Xi(\mathbf{q})} \sum_{k \in \Omega_j} \left(w_j(\mathbf{q})\psi_k(\mathbf{q})\right) \dot{\varphi}_k. \tag{28}$$

### 3.3 Nonlinear system of collocation equations

We shall now tackle the collocation of the transformed Navier-Stokes equations. The closed domain $\overline{\Omega}_T$ is discretised into the set[2] $Q_I$ of $N_I$ interior nodes $Q_I = \{\mathbf{q}_1,\ldots,\mathbf{q}_{N_I}\}$, plus the set $Q_B$ of $N_B$ nodes along the boundary, $Q_B = \{\mathbf{q}_{1+N_I},\ldots,\mathbf{q}_N = \mathbf{q}_{N_B+N_I}\}$. In turn, the set of boundary nodes $Q_B$ is ordered as follows: first, the set $Q_{FAR}$ of far-field nodes ($\xi = \ell$); then, the set $Q_{CYL}$ of nodes on the cylinder, or near-field nodes; and lastly the set $Q_{AXIS}$ of nodes on the $x$-axis (with either $\varphi = 0$ or $\varphi = \pi$). The corresponding number of nodes is #$FAR$, #$CYL$ and #$AXIS$ respectively, thus $N_B =$ #$FAR +$ #$CYL +$ #$AXIS$.

The sought-for unkowns are the nodal values of the the RBF-PU interpolants, which we choose to assemble into a *nodal system vector*

$$\vec{X} := \left(\dot{\xi}_1,\ldots,\dot{\xi}_N, \dot{\varphi}_1,\ldots,\dot{\varphi}_N, p_1,\ldots,p_N\right)^T, \tag{29}$$

where $\dot{\varphi}_j$ is the RBF-PU interpolant for $\dot{\varphi}$ evaluated at node $\mathbf{q}_j$, etc.

According to (25)—(27), collocation of the BCs in Figure 3 yields the following matrix $[BC]$ of size $\#BC \times 3N$, where $\#BC = 3(\#FAR + \#AXIS) + 2\#CYL$:

$$\begin{pmatrix} [1]_{FAR} & 0 & 0 \\ 0 & [1]_{FAR} & 0 \\ 0 & 0 & [1]_{FAR} \\ [1]_{CYL} & 0 & 0 \\ 0 & [1]_{CYL} & 0 \\ [\partial_\varphi]_{AXIS} & 0 & 0 \\ 0 & [1]_{AXIS} & 0 \\ 0 & 0 & [\partial_\varphi]_{AXIS} \end{pmatrix} \vec{X} =: [BC]\vec{X} = \vec{g} := \begin{pmatrix} (\cos\varphi_1,\ldots,\cos\varphi_{\#FAR})^T \\ (-\sin\varphi_1,\ldots,-\sin\varphi_{\#FAR})^T \\ 0 \\ \vdots \\ 0 \end{pmatrix}. \tag{30}$$

---

[2]To simplify set notation, we shall write that $j \in Q$ if $\mathbf{q}_j \in Q$, where $Q$ is any of the defined node sets.

The collocation of the nonlinear operators $\mathscr{W}_1$ and $\mathscr{W}_2$ on the interior nodes gives rise to the $2N_I$ equations nonlinear in $\vec{X}$:

$$\mathscr{W}_1(\mathbf{q}_1,\vec{X}) = 0,\ldots,\mathscr{W}_1(\mathbf{q}_{N_I},\vec{X}) = 0, \mathscr{W}_2(\mathbf{q}_1,\vec{X}) = 0,\ldots,\mathscr{W}_2(\mathbf{q}_{N_I},\vec{X}) = 0, \quad (31)$$

where $\mathscr{W}(\mathbf{q}',\vec{X}') \in \mathbb{R}$ means that the three RBF-PU interpolants and those of their derivatives are computed at nodal values $\vec{X}'$, combined according to the operator $\mathscr{W}$, and collocated at $\mathbf{q}'$.

On the other hand, $\mathscr{W}_3$ is a linear operator. Note that, by (8),

$$\mathscr{W}_3(\mathbf{q}_k,\vec{X}) = \dot{\xi}_k + (\ell - \xi_k)\sum_{j=1}^{N}[\partial_\xi]_{kj}\dot{\xi}_j + \sum_{j=1}^{N}[\partial_\varphi]_{kj}\dot{\varphi}_j. \quad (32)$$

There is no BC for the pressure on the cylinder (for $\vec{p}$ on $Q_{CYL}$). In order to have as many equations as unknowns in $\vec{X}$, we enforce the linear PDE $\mathscr{W}_3$ (representing the incompressibility of the fluid) also on the surface of the cylinder. In this way, the collocation of $\mathscr{W}_3$ gives rise to the following block of $N_I + \#CYL$ linear equations in $\vec{X}$:

$$\left(\begin{array}{c|c|c} [1]_{Q_I} + diag_{Q_I}(\ell - \xi)[\partial_\xi]_{Q_I} & [\partial_\varphi]_{Q_I} & [0] \\ \\ [1]_{Q_{CYL}} + diag_{Q_{CYL}}(\ell - \xi)[\partial_\xi]_{Q_{CYL}} & [\partial_\varphi]_{Q_{CYL}} & [0] \end{array}\right)\vec{X} =: [\mathscr{W}_3]\vec{X} = \vec{0}, \quad (33)$$

where the $[0]$ are padding matrices of zeros of the appropriate size, and the diagonal matrix $diag_Q(f(\xi,\varphi))$ is defined as follows. Let $Q = \{\mathbf{q}'_1,\ldots,\mathbf{q}'_{\#Q}\}$, then

$$diag_Q\big(f(\xi,\varphi)\big) := \begin{pmatrix} f(\xi'_1,\varphi'_1) & 0 & \ldots & 0 \\ 0 & f(\xi'_2,\varphi'_2) & \ldots & 0 \\ \vdots & \vdots & \ddots & \vdots \\ 0 & \ldots & 0 & f(\xi'_{\#Q},\varphi'_{\#Q}) \end{pmatrix}. \quad (34)$$

Summing up, $\vec{X}$ solves the square system mixing linear and nonlinear equations

$$\begin{cases} \left.\begin{array}{l} [BC]\vec{X} = \vec{g} \\ [\mathscr{W}_3]\vec{X} = (0,\ldots,0)^T \end{array}\right\} \text{ linear block } [LIN]\vec{X} = (\vec{g}^T,0,\ldots,0)^T =: \vec{g}_{LIN} \\ \mathscr{W}_1(\mathbf{q}_1,\vec{X}) = 0 \\ \vdots \\ \mathscr{W}_1(\mathbf{q}_{N_I},\vec{X}) = 0 \\ \mathscr{W}_2(\mathbf{q}_1,\vec{X}) = 0 \\ \vdots \\ \mathscr{W}_2(\mathbf{q}_{N_I},\vec{X}) = 0 \end{cases} \quad (35)$$

**Elimination of the BCs.** It is advantageous to eliminate the $N_I + 3N_B$ linear equations before solving the nonlinear system, thus shrinking its size. The optimally stable way of doing so is by means of the QR decomposition of $[LIN]$—see [28, section 15.2]:

$$[LIN]^T \Pi = [O_1 O_2] \begin{bmatrix} R \\ 0 \end{bmatrix}, \tag{36}$$

where $\Pi$ is a permutation matrix, R is upper triangular, and $O_1 \in \mathbb{R}^{3N \times (N_I + 3N_B)}$ and $O_2 \in \mathbb{R}^{3N \times 2(N-N_B)}$ are made up of orthogonal columns.[3] Then, $\vec{X}$ can be expressed in terms of a fixed vector and a smaller vector $\vec{Y}$ (which remains to be found) as

$$\vec{X} = O_1 R^{-T} \Pi^T \vec{g}_{LIN} + O_2 \vec{Y}. \tag{37}$$

RBF-PU interpolants depend now on $\vec{Y}$ via (37). The smaller, purely nonlinear final set of equations $\{E_i(\vec{Y})\}_{i=1}^{i=2N_I}$ for the unknowns $\vec{Y}$ is

$$\begin{cases} E_1(\vec{Y}) = \mathscr{W}_1(\mathbf{q}_1, \vec{Y}) = 0 \\ \vdots \\ E_{N_I}(\vec{Y}) = \mathscr{W}_1(\mathbf{q}_{N_I}, \vec{Y}) = 0 \\ E_{1+N_I}(\vec{Y}) = \mathscr{W}_2(\mathbf{q}_1, \vec{Y}) = 0 \\ \vdots \\ E_{2N_I}(\vec{Y}) = \mathscr{W}_2(\mathbf{q}_{N_I}, \vec{Y}) = 0. \end{cases} \tag{38}$$

## 3.4 The trust region algorithm and the RBF-PU Jacobian

In order to solve for $\vec{Y}$, we apply the trust region algorithm, which transforms a rootfinding problem (the root being the vector solution of the system) into a minimization problem for the (half) sum-of-squares residual, or *merit function*, $\mu(\vec{Y}') := (1/2) \sum_{i=1}^{2N_I} E_i^2(\vec{Y}')$. Because $\mu$ can be highly nonconvex, it is critical both for convergence and for speed that the analytic Jacobian be available. The entry $J_{ij}$ of the Jacobian at $\vec{Y} = \vec{Y}'$ is

$$J_{ij}(\vec{Y}') = \frac{\partial E_i}{\partial Y_j}(\vec{Y}'), \qquad 1 \leq i, j \leq 2N_I. \tag{39}$$

For clarity, we shall discuss the Jacobian with respect to $\vec{X}$ first. The derivatives of the collocated nonlinear equation with respect to the RBF-PU interpolant nodal values in $\vec{X}$ are calculated via the Fréchet derivatives of the nonlinear operators $\mathscr{W}_1$ and $\mathscr{W}_2$ with respect to linear operators $\mathscr{L}$. Then, the latter are replaced by those of the RBF-PU interpolant at the collocation node. Formally, the sum includes all the derivatives up to second order, but many of them are zero. The precise dependencies are:

$$\mathscr{W}_1 = \mathscr{W}_1\left(\dot{\xi}, \dot{\varphi}, \partial_\xi \dot{\xi}, \partial_\varphi \dot{\xi}, \partial_\varphi \dot{\varphi}, \partial_\xi p, \partial_{\xi\xi}^2 \dot{\xi}, \partial_{\varphi\varphi}^2 \dot{\xi}\right), \tag{40}$$

$$\mathscr{W}_2 = \mathscr{W}_2\left(\dot{\xi}, \dot{\varphi}, \partial_\varphi \dot{\xi}, \partial_\varphi \dot{\varphi}, \partial_\xi \dot{\varphi}, \partial_\varphi p, \partial_{\xi\xi}^2 \dot{\varphi}, \partial_{\varphi\varphi}^2 \dot{\varphi}\right). \tag{41}$$

---

[3] Since $Q$ has is repeatedly used to denote sets, we call orthogonal matrices $O$.

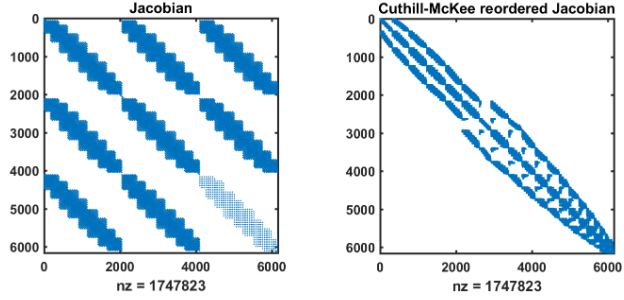

Figure 4: Sparsity pattern of the Jacobian matrix (53) corresponding to Figure 5. Instead of the QR procedure in Section 3.3, the Dirichlet BCs were directly eliminated (left). On the right, same matrix after applying Cuthill-McKee permutations.

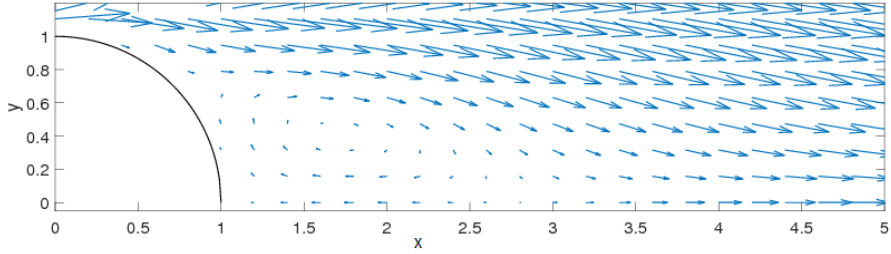

Figure 5: Velocity vectors in the wake region past a circular cylinder at $Re = 20$.

The nonzero Fréchet derivatives of $\mathscr{W}_1$ and $\mathscr{W}_2$ are, from (6)-(7):

$$\frac{\partial \mathscr{W}_1}{\partial(\partial^2_{\xi\xi}\dot{\xi})} = \frac{\partial \mathscr{W}_2}{\partial(\partial^2_{\xi\xi}\dot{\varphi})} = -\frac{(\ell - \xi)^3}{\ell}, \qquad \frac{\partial \mathscr{W}_1}{\partial(\partial^2_{\varphi\varphi}\dot{\xi})} = \frac{\partial \mathscr{W}_2}{\partial(\partial^2_{\varphi\varphi}\dot{\varphi})} = -1 + \frac{\xi}{\ell}, \quad (42)$$

$$\frac{\partial \mathscr{W}_1}{\partial \dot{\xi}} = \frac{Re}{2}(\ell - \xi)\partial_\xi \dot{\xi} + 1 - \frac{\xi}{\ell}, \qquad \frac{\partial \mathscr{W}_2}{\partial \dot{\xi}} = \frac{Re}{2}\left((\ell - \xi)\partial_\xi \dot{\varphi} + \dot{\varphi}\right), \quad (43)$$

$$\frac{\partial \mathscr{W}_1}{\partial \dot{\varphi}} = \frac{Re}{2}\left(\partial_\varphi \dot{\xi} - 2\dot{\varphi}\right), \qquad \frac{\partial \mathscr{W}_2}{\partial \dot{\varphi}} = \frac{Re}{2}\left(\dot{\xi} + \partial_\varphi \dot{\varphi}\right) + 1 - \frac{\xi}{\ell}, \quad (44)$$

$$\frac{\partial \mathscr{W}_1}{\partial(\partial_\xi \dot{\xi})} = \frac{Re}{2}(\ell - \xi)\dot{\xi} + \frac{(\ell - \xi)^2}{\ell} = \frac{\partial \mathscr{W}_2}{\partial(\partial_\xi \dot{\varphi})}, \quad (45)$$

$$\frac{\partial \mathscr{W}_1}{\partial(\partial_\varphi \dot{\xi})} = \frac{Re}{2}\dot{\varphi}, \qquad \frac{\partial \mathscr{W}_2}{\partial(\partial_\varphi \dot{\xi})} = -2\left(1 - \frac{\xi}{\ell}\right), \quad (46)$$

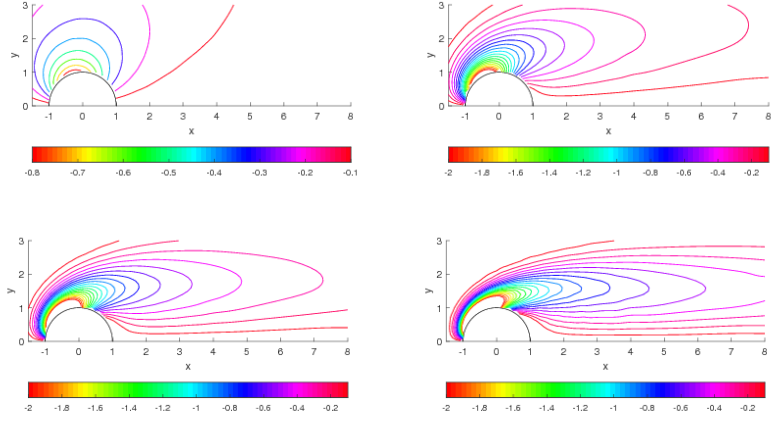

Figure 6: Vorticity contours of the steady flow past a circular cylinder at various Reynolds numbers: $Re = 1$ (top left; note the different scale), $Re = 10$ (top right) $Re = 20$ (bottom left) and $Re = 40$ (bottom right). The contour levels are -2:0.1:-0.1.

$$\frac{\partial \mathscr{W}_1}{\partial(\partial_\varphi \dot\varphi)} = 2\left(1 - \frac{\xi}{\ell}\right), \qquad \frac{\partial \mathscr{W}_2}{\partial(\partial_\varphi \dot\varphi)} = \frac{Re}{2}\dot\varphi, \tag{47}$$

$$\frac{\partial \mathscr{W}_1}{\partial(\partial_\xi p)} = \frac{Re}{2}(\ell - \xi), \qquad \frac{\partial \mathscr{W}_2}{\partial(\partial_\varphi p)} = \frac{Re}{2}. \tag{48}$$

If $E_i$ is the collocation of $\mathscr{W}_1$ on $\mathbf{q}_k$, it holds, by the chain rule and from (29):

- If $1 \leq j \leq N$ (so that $X_j = \dot\xi_j$), then

$$\frac{\partial E_i}{\partial X_j} = \left(\frac{\partial \mathscr{W}_1}{\partial \dot\xi}\right)(\mathbf{q}_k, \vec{X}')[1]_{kj} + \left(\frac{\partial \mathscr{W}_1}{\partial(\partial_\xi \dot\xi)}\right)(\mathbf{q}_k, \vec{X}')[\partial_\xi]_{kj} + \left(\frac{\partial \mathscr{W}_1}{\partial(\partial_\varphi \dot\xi)}\right)(\mathbf{q}_k, \vec{X}')[\partial_\varphi]_{kj}$$
$$+ \left(\frac{\partial \mathscr{W}_1}{\partial(\partial^2_{\xi\xi} \dot\xi)}\right)(\mathbf{q}_k, \vec{X}')[\partial^2_{\xi\xi}]_{kj} + \left(\frac{\partial \mathscr{W}_1}{\partial(\partial^2_{\varphi\varphi} \dot\xi)}\right)(\mathbf{q}_k, \vec{X}')[\partial^2_{\varphi\varphi}]_{kj}. \tag{49}$$

- If $N+1 \leq j \leq 2N$ (so that $X_j = \dot\varphi_j$), then

$$\frac{\partial E_i}{\partial X_j} = \left(\frac{\partial \mathscr{W}_1}{\partial(\partial_\varphi \dot\varphi)}\right)(\mathbf{q}_k, \vec{X}')[\partial_\varphi]_{kj}. \tag{50}$$

- Finally, if $2N+1 \leq j \leq 3N$ (so that $X_j = p_j$), then

$$\frac{\partial E_i}{\partial X_j} = \left(\frac{\partial \mathscr{W}_1}{\partial(\partial_\xi p)}\right)(\mathbf{q}_k, \vec{X}')[\partial_\xi]_{kj}. \tag{51}$$

Repeating the procedure for $\mathcal{W}_2$ yields the lower block of $J(\vec{X})$. Once the entries $\partial E_i / \partial X_j$ are available, the $\partial E_i / \partial Y_j$ are computed by applying the chain rule to (37).

The Jacobian can be compactly written as the square matrix (see [3] for details):

$$J(\vec{Y}) = \sum_{k \in NZFD} \begin{pmatrix} diag_{Q_l}\left(\frac{\partial \mathcal{W}_1}{\partial \mathcal{L}_k \xi}(\vec{Y})\right) \cdot [\mathcal{L}_k] & diag_{Q_l}\left(\frac{\partial \mathcal{W}_1}{\partial \mathcal{L}_k \phi}(\vec{Y})\right) \cdot [\mathcal{L}_k] & diag_{Q_l}\left(\frac{\partial \mathcal{W}_1}{\partial \mathcal{L}_k p}(\vec{Y})\right) \cdot [\mathcal{L}_k] \\ diag_{Q_l}\left(\frac{\partial \mathcal{W}_2}{\partial \mathcal{L}_k \xi}(\vec{Y})\right) \cdot [\mathcal{L}_k] & diag_{Q_l}\left(\frac{\partial \mathcal{W}_2}{\partial \mathcal{L}_k \phi}(\vec{Y})\right) \cdot [\mathcal{L}_k] & diag_{Q_l}\left(\frac{\partial \mathcal{W}_2}{\partial \mathcal{L}_k p}(\vec{Y})\right) \cdot [\mathcal{L}_k] \end{pmatrix} O_2 .$$

(NZFD= nonzero Fréchet derivatives in $\mathcal{W}_1$ and $\mathcal{W}_2$)

$$=: \hat{J}(\vec{Y}) O_2. \tag{52}$$

Aside from $O_2$, all the submatrices above are either diagonal or sparse matrices. Upon multiplication by $O_2$, however, that sparse pattern for $J(\vec{Y})$ is lost.

**The dogleg approximation.** The trust region algorithm [28] is an heuristic framework for finding $\mu(\vec{Y}_{root}) = 0$, starting from an initial guess $\vec{Y}_0$ and iterating such that (hopefully) $\lim_{k \to \infty} \vec{Y}_k = \vec{Y}_{root}$. The iterates are constructed according to the Hessian of the merit function, $H(\vec{Y}_k)$. The simplest implementation is the *dogleg method*, in which $H(\vec{Y}_k) \approx J^T(\vec{Y}_k) J(\vec{Y}_k)$. It is particularly convenient[4] in the RBF context since steps can be constructed with matrix-vector products and linear systems involving the matrix $J(\vec{Y}_k)$ only—hence at half of the condition number order of $J^T(\vec{Y}_k) J(\vec{Y}_k)$.

By choosing a proper initial guess (as explained in Section 4), all the experiments carried out in this work converged to what distinctly seems a unique root of the system of nonlinear RBF-PU collocation equations. In tougher scenarios, a more robust variant of the trust region algorithm than the dogleg method may be required, which uses a better approximation to the Hessian. The reader is directed to [3] for details.

**An alternative.** The reason why $J(\vec{Y})$ is not sparse is the multiplication by $O_2$, which stems from the optimal elimination of linear equations from the nonlinear system (37).

An alternative approach is eliminating *only the Dirichlet BCs* directly from (30). This yields a vector of unknowns $\vec{X}'$, smaller than $\vec{X}$ but larger than $\vec{Y}$, since the collocation equations for the Neumann BCs as well as for $\mathcal{W}_3$ are retained in the system. Let $[LIN']\vec{X}'$ be the remaining non-Dirichlet linear equations. The Jacobian for $\vec{X}'$ reads

$$J(\vec{X}') = \left( \frac{[LIN']}{\hat{J}(\vec{X}')} \right) \in \mathbb{R}^{(\#\vec{X}') \times (\#\vec{X}')}, \tag{53}$$

where $\hat{J}(\vec{X}')$ is analogous to the matrix block in (52). While larger than $J(\vec{Y})$, $J(\vec{X}')$ is sparse, as can be seen in Figure 4. In our experiments, that Jacobian was marginally worse conditioned than, and performed similarly to, $J(\vec{Y})$. The triply banded structure can be permuted into a diagonal sparsity profile with minimal bandwidth using the sparse reverse Cuthill-McKee algorithm (implemented in Matlab with the `symrcm` command), even though $J(\vec{X}')$ is not symmetric—see Figure 4.

---

[4]In a Matlab implementation such as ours, `fsolve` must be fed the system (38) (easy to code from (35) and the mapping (37)); the Jacobian (52); and a guess $\vec{Y}_0$ of the solution to kick off the trust-region iterations. (And the trust region algorithm must be set as solver.)

# 4  Steady flow past a circular cylinder

We restrict ourselves to $Re \leq 40$, before the onset of instability at about $Re \approx 47$. For computations, most authors truncate the domain at a "large enough" radius and apply some unphysical numerical BC at the artificial interface. It is well documented that this can affect the solution near the boundary of the obstacle [27]. We instead consider a semi-infinite domain (using a reflecting BC on $y = 0$) and impose the actual BCs at infinity.

In order to compute the steady flow, practitioners typically solve the time-dependent Navier-Stokes equations to equilibrium. In that way, nonlinearities can be circumvented by an explicit timestepping scheme, whose stability may nonetheless dictate small timesteps—thus leading to lengthy overall simulations. Instead, we directly solve the nonlinear elliptic equations which model the equilibrium flow, and benefit from very precise Jacobians, thanks to the accuracy of RBF-PU in reproducing derivatives. In order to define a good initial guess for the trust-region algorithm, the numerical solution for a lower $Re$ is a natural and robust initial choice. For instance, in order to solve the flow at $Re = 40$, we followed the sequence of solves

$$Re = 1 \text{ (initial guess identically zero)} \quad \rightarrow \quad Re = 20 \quad \rightarrow \quad Re = 40.$$

The RBF-PU solution was captured in at most nine iterations in all cases.[5] (On the other hand, leaping from $Re = 1$ to $Re = 40$ would not converge.)

After some preliminary tests, we concluded that the variations induced by the stretching scale $\ell$ in (3) are small, and the medium scale $\ell = 2$ was chosen for our simulations. They were performed in a rectangular transformed domain $[0,2] \times [0,\pi]$, discretised into a grid of RBF nodes a distance $h$ apart, with $0.05 \leq h \leq 0.10$ across the simulations. The largest number of RBFs was 2520 (with $h = 0.05$), with three times as many nodal values (for $\dot{\xi}, \dot{\phi}$, and $p$). The cover consisted of gridded circular patches with radii $\rho_i = 0.25$ and centres separated by one radius. The Jacobian matrix sparsity ratio increases with the patches' radii (see Figure 4 for illustration). In all cases, the condition number of the Jacobian matrix was about $10^9$, where the RBF was the inverse multiquadric, $\phi(r) = (1 + \varepsilon^2 r^2)^{-1/2}$ with the shape parameter set to $\varepsilon = 2$.

*Remark.* Picking a "near optimal" shape parameter is both a challenging and an advantageous feature of RBF methods, but a topic that we have chosen to avoid in this paper. Thus, we have used the same, unoptimised value $\varepsilon = 2$ in all of our simulations.

Translating the solution back into the physical domain, we obtained flow patterns such as that one shown in Figure 5. The bubble of recirculation immediately past the cylinder is well resolved; its size increases with *Re* as expected [12].

---

[5] It must be stressed that convergence to the root of the nonlinear system of collocation equations might well not take place at all (particularly at higher $Re$) with a less robust nonlinear solver than the trust region algorithm, such as Levenberg-Marquardt's or Newton's method. In fact, convergence would take plenty of iterations, or get stuck in a local minimum even with the trust region algorithm, but using expensive finite-differences Jacobians instead of the RBF-PU semianalytical one (52). (As happens by default in Matlab's `fsolve` when no Jacobian is provided.) The reason is that the inaccuracy of the finite difference Jacobian is compounded by the typically high condition numbers in RBF simulations. On the other hand, we did not observe any such difficulty in any of the experiments in this paper, both in this section and in the next one.

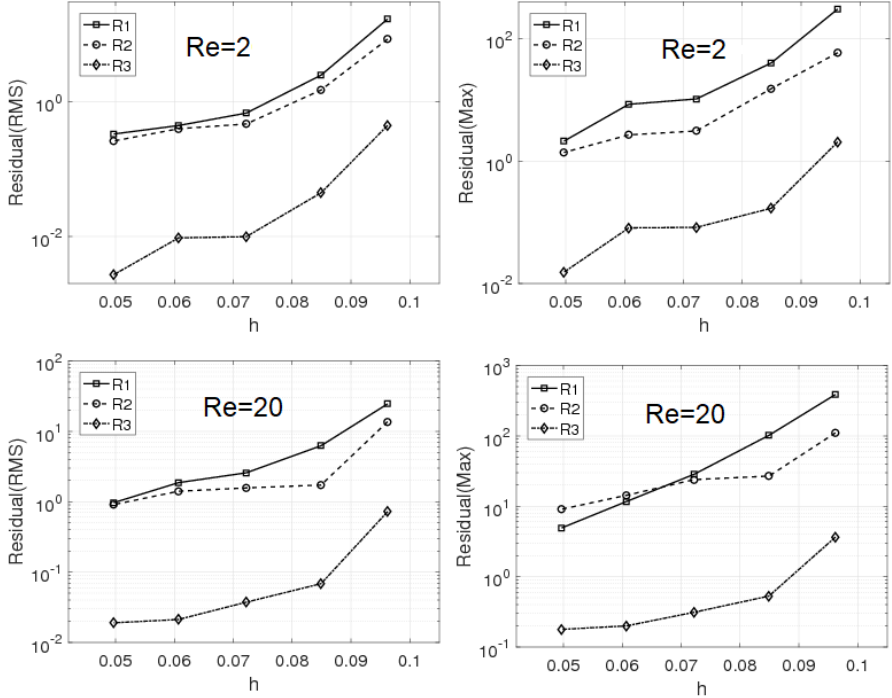

Figure 7: The RMS and Max residual of RBF-PU solution of circular cylinder at $Re = 2$ and $Re = 20$ as a function of nodal distance $h$. $R_i$ are the residuals to $\mathscr{W}_i$, for $i = 1, 2, 3$.

Figure 6 shows the vorticity contours of the numerical flow. The vorticity function, $\omega = |\nabla \times \mathbf{u}|$, reads in the transformed domain

$$\omega = \tfrac{\ell-\xi}{\ell}\left((\ell-\xi)\partial_\xi \dot\varphi - \partial_\varphi \dot\xi + \dot\varphi\right), \quad 0 \le \xi \le \ell, 0 \le \varphi \le \pi. \tag{54}$$

The quality of the RBF-PU solution can be assessed through the root mean square (RMS) value of the residual.[6] Let $R_1$, $R_2$, $R_3$ denote the pointwise residual of RBF-PU solutions to the equations (5). The RMS and Max metrics with respect to $h$ are displayed in Figure 7; both metrics of the three residuals decay as resolution (in the transformed domain) increases. Due to the fact that—with smooth solutions—pointwise error and residual are correlated [8], convergence of the error with $h$ is to be expected.

We continue our study by comparing with numerical results in the literature.

**Comparisons and drag coefficient.** Let $p|_\Gamma = p/(\tfrac{1}{2}U_\infty)$ and $\omega|_\Gamma$ denote respectively the pressure and the vorticity on the cylinder surface. The profiles of both of them are shown in Figure 8. (They bear a very close resemblance to figure 1 in reference [34].)

---

[6]It was sampled over a set of evaluation nodes $[-2:0.2:8] \times [0:0.2:5]$ in the physical domain.

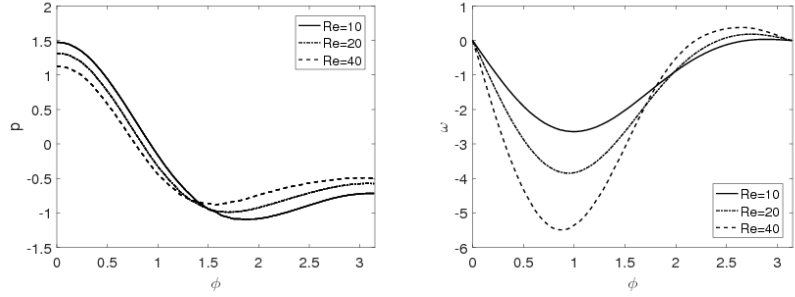

Figure 8: The pressure (left) and the vorticity (right) on the cylinder surface. $\phi = \pi - \varphi$, so that $\phi = 0$ is the front stagnation point, of physical coordinates $(-R = -1, 0)$.

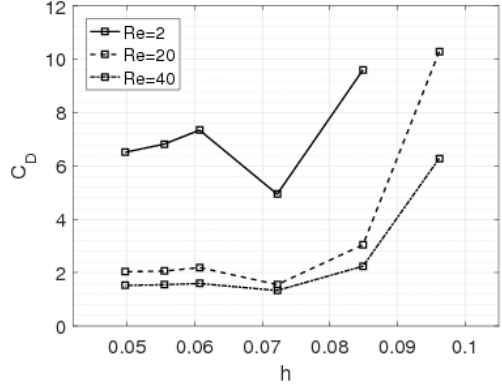

Figure 9: Behaviour of the drag coefficient of the circular cylinder as a function of $h$.

Let $\tau = (\tau_x, \tau_y)$ be the unit tangent vector on the surface and $\mathbf{n} = (n_x, n_y) = (\tau_y, -\tau_x)$ the outward normal. The drag coefficient of the immersed cross section is defined as

$$C_D = \underbrace{\int_0^\pi p|_\Gamma \, n_x \, d\varphi}_{C_p} - \underbrace{\frac{4}{Re} \int_0^\pi \omega|_\Gamma \, \tau_x \, d\varphi}_{C_\omega}. \tag{55}$$

Above, the contribution of $C_p$ to the drag is due to pressure, while $C_\omega$ is due to the viscosity. For the circular cylinder,

$$C_D = \int_0^\pi p|_{\xi=0} \cos\varphi \, d\varphi - \frac{4}{Re} \int_0^\pi \omega|_{\xi=0} \sin\varphi \, d\varphi. \tag{56}$$

According to Table 1, the RBF-PU results are in good agreement with simulations[7]

---

[7] $L$ is the distance from the rear of the cylinder to the end of the recirculation region. It can be measured by finding the point where the horizonal velocity $u$ along the $x$-axis vanishes. The length is reported normalised with respect to the cylinder diameter.

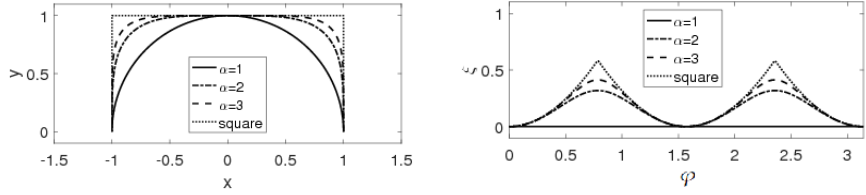

Figure 10: Profiles of the rounded square cylinder for several values of $\alpha$, and its shape in the transformed domain under the nonlinear change of variables. Even with the cusps, the transformed domain complies with an interior cone condition [37].

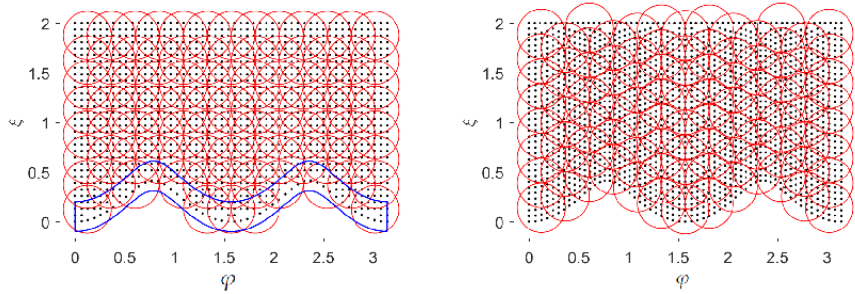

Figure 11: Gridded (left) and body-fitted (right) discretisations of RBF-PU nodes and patches of the computational domain for $\alpha = 3$. Note the gaps in the region near the boundary (inside the blue band)—critical for reproducing the fluid-body interaction.

in [11, 22, 34, 5, 13, 23] and experimental data reported in [7]. Note that all other numerical results used truncated computational domains.

The convergence of the drag coefficient with $h$ is depicted in Figure 9. It can be seen that $C_D$ becomes essentially independent of the discretisation rather fast.

## 5 Steady flow past rounded or sharp square cylinders

### 5.1 Discretisation

An equation for rounded square shapes with various degrees of smoothness is

$$x^{2\alpha} + y^{2\alpha} = 1, \quad \alpha = 1, 2, 3, \ldots \qquad (57)$$

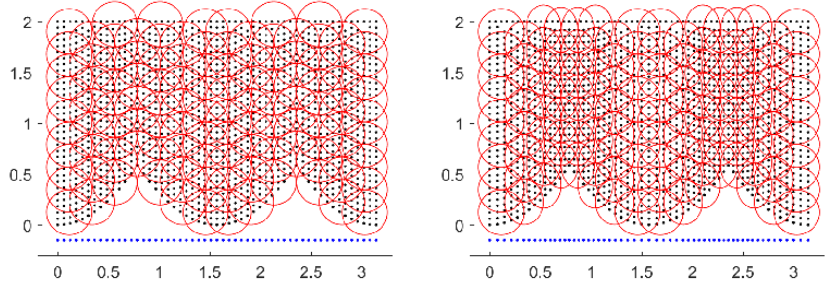

Figure 12: Two body-fitted RBF-PU discretisations of the square cylinder. The one on the right clusters nodes around the cusps according to the distribution of blue points depicted along the horizonal ($\varphi$) axis. Patches are squeezed horizontally too, becoming ellipses, so that the number of nodes per patch is roughly constant.

Table 1: Comparison of characteristic dimensions of the flow past a circular cylinder at $Re = \{20, 40\}$. The pressure, viscous and drag coefficients are $C_p$, $C_\omega$ and $C_D$, respectively. $L$ is the wake length measused from the rear stagnation point on the surface, of physical coordinates $(1,0)$ (in our simulations, the circle centre is the origin). $(a,b)$ is the location of the eddy centre from $(1,0)$. All lengths are measured in diameters. Data from [7] are experimental. Dashes correspond to unavailable results.

|  | Re | $C_\omega$ | $C_p$ | $C_D$ | $L$ | $(a,b)$ |
| --- | --- | --- | --- | --- | --- | --- |
| *Ståhlberg et al.* (2006) [34] | 20 | 0.82 | 1.23 | 2.05 | 0.90 | – |
|  | 40 | 0.54 | 0.99 | 1.53 | 2.13 | – |
| *Linnick & Fasel* (2005) [22] | 20 | – | – | 2.16 | 0.93 | (0.36, 0.43) |
|  | 40 | – | – | 1.61 | 2.23 | (0.71, 0.59) |
| *Liu* (2017) [23] | 20 | – | – | 2.09 | 0.98 | (0.37, 0.43) |
|  | 40 | – | – | 1.59 | 2.40 | (0.74, 0.60) |
| *Fornberg* (1980) [11] | 20 | – | – | 2.00 | 0.91 | – |
|  | 40 | – | – | 1.50 | 2.24 | – |
| *Gautier et al.* (2013) [13] | 40 | – | – | 1.49 | 2.24 | (0.71, 0.59) |
| *Bouchon et al.* (2012) [5] | 40 | – | – | 1.50 | 2.26 | (0.71, 0.60) |
| *Coutanceau & Bouard* (1977) [7] | 40 | – | – | – | 2.13 | (0.76, 0.59) |
| *current study* | 20 | 0.80 | 1.23 | 2.03 | 0.91 | (0.36, 0.42) |
|  | 40 | 0.50 | 1.02 | 1.52 | 2.17 | (0.72, 0.60) |

After applying (3), the rounded square in the transformed domain is described as

$$\xi = \ell \left[ 1 - (\cos^{2\alpha} \varphi + \sin^{2\alpha} \varphi)^{\frac{1}{2\alpha}} \right]. \tag{58}$$

The rounded square in the physical and transformed domains is depicted for several values of $\alpha$ in Figure 10 ($\alpha = 1$ is the circular cylinder). In a first attempt at discretising the transformed domain, we started with the same grid of RBF-PU nodes and patches as we did with the circular cylinder. However, with the rounded square cylinder, this discretisation does not sample the region in the transformed domain close to the cylinder surface adequately, owing to the resulting gaps. This leads to an unacceptable loss of accuracy of the RBF-PU near fields. This issue can be fixed by adjusting the position of patches along the transformed boundary according to (58), as well (see Figure 11). While the straightforward discretisation of nodes and patches (left) does not evenly sample the interesting region (highlighted in blue), the adaptive one (right) helps capture it much better. (Note that the gaps have been pushed towards $\xi = 2$, i.e. infinitely far into the physical domain.)

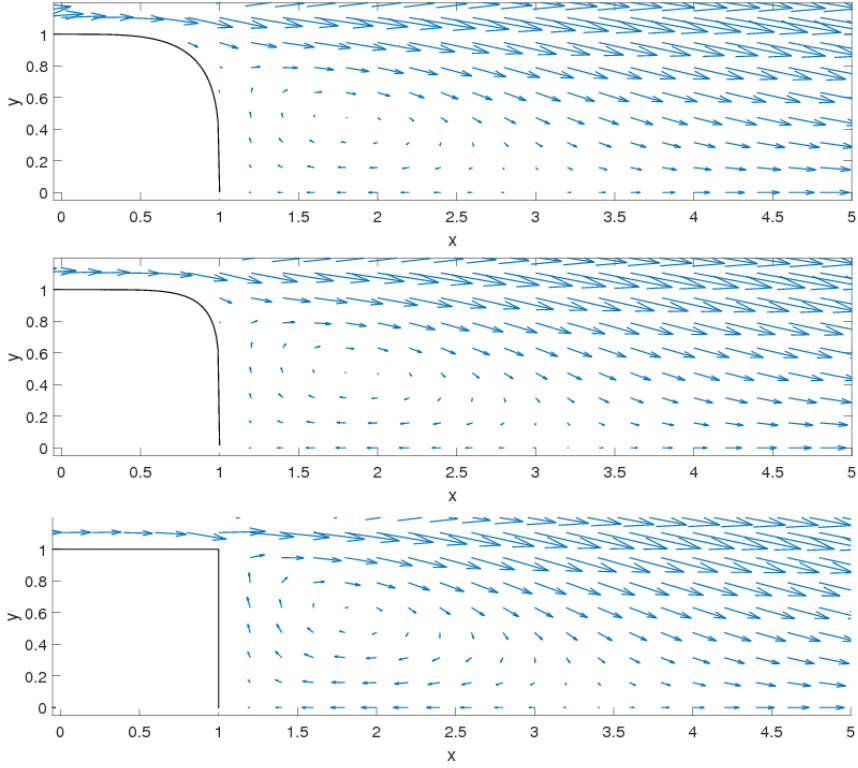

Figure 13: Flow pattern in the wake for $Re = 20$ and $\alpha = \{2, 3, \infty\}$ (from top to bottom).

The limit $\alpha \to \infty$ of the transformation (58) is the perfectly square cylinder (see Figure 10). More practically, the latter can be mapped onto the transformed domain by

$$\xi = \begin{cases} \ell(1 - \cos\varphi) & \text{if } 0 \leq \varphi < \pi/4, \\ \ell(1 - \sin\varphi) & \text{if } \pi/4 \leq \varphi < 3\pi/4, \\ \ell(1 - \cos(\pi - \varphi)) & \text{if } 3\pi/4 \leq \varphi \leq \pi. \end{cases} \qquad (59)$$

The nodes are not evenly located along the $\varphi$-axis in the transformed domain, but

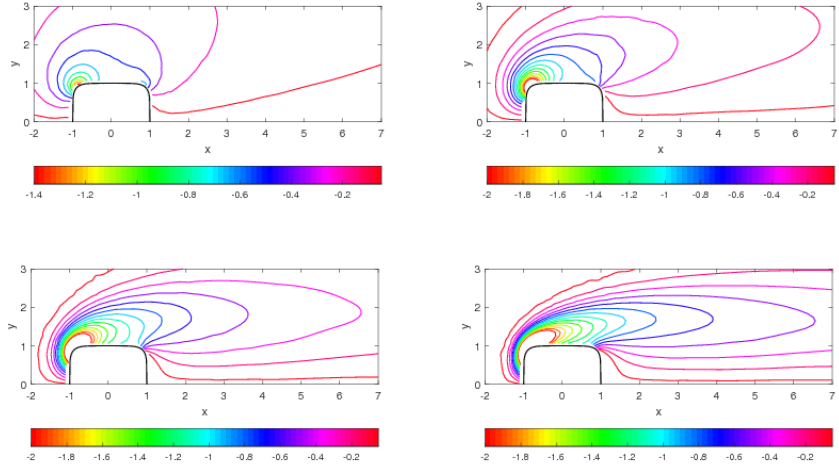

Figure 14: The vorticity contours for the steady flow past rounded square cylinder $\alpha = 3$ for different values of Reynolds number. $Re = 2$ (top left), $Re = 10$ (top right) $Re = 20$ (below left) and $Re = 40$ (below right).

clustered near the cusps according to the following ad-hoc formula:

$$\varphi = \left\{ \begin{array}{c} \pi/4 \\ 3\pi/4 \end{array} \right\} + \lambda \sinh(\eta_j), \text{ where } \eta_j \in \left[ \sinh^{-1}\left(\frac{-\pi}{4\lambda}\right), \sinh^{-1}\left(\frac{\pi}{4\lambda}\right) \right]. \quad (60)$$

The $\{\eta_j\}$ are equispaced and $\lambda$ is a parameter that determines the refinement at corner points. The patch centres are defined with a similar pattern as for the node distribution. Moreover, they are "squeezed" into ellipses—so that all patches contain approximately as many nodes—as discussed in [30]. See Figure 12.

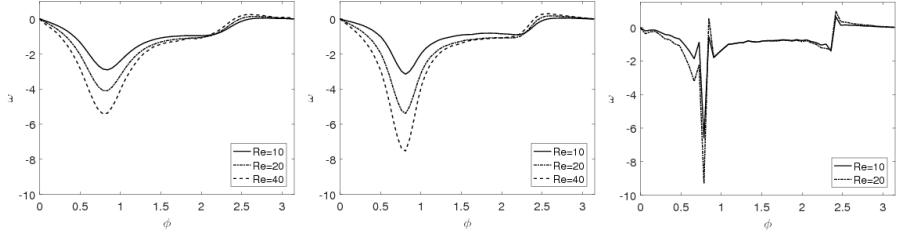

Figure 15: Vorticity profile along the obstacle for $\alpha = \{2, 3, \infty\}$ (from left to right). Note $\phi = \pi - \varphi$, hence $\omega|_\Gamma(\phi)$ is shown from the front to the rear stagnation points.

Table 2: Characteristic flow values (*L*, *a* and *b* measured in square side lengths).

|  | Re | $C_\omega$ | $C_p$ | $C_D$ | L | (a,b) |
|---|---|---|---|---|---|---|
| **rounded cylinder ($\alpha = 2$)** | | | | | | |
| *current study* | 10 | 0.97 | 1.37 | 2.34 | 0.39 | (0.18, 0.24) |
|  | 20 | 0.68 | 1.10 | 1.78 | 0.99 | (0.41, 0.43) |
|  | 30 | 0.48 | 0.94 | 1.42 | 1.68 | (0.57, 0.60) |
|  | 40 | 0.34 | 0.84 | 1.19 | 1.98 | (0.68, 0.60) |
| **rounded cylinder ($\alpha = 3$)** | | | | | | |
| *current study* | 10 | 0.86 | 1.25 | 2.10 | 0.38 | (0.22, 0.28) |
|  | 20 | 0.64 | 1.17 | 1.81 | 1.14 | (0.47, 0.46) |
|  | 30 | 0.47 | 1.00 | 1.47 | 2.03 | (0.61, 0.56) |
|  | 40 | 0.38 | 0.94 | 1.33 | 2.59 | (0.70, 0.63) |
| **square cylinder ($\alpha = \infty$)** | | | | | | |
| *Saha* (2013) [31] | 10 | – | – | 3.07 | 0.64 | (–, 0.49) |
|  | 20 | – | – | 2.21 | 1.35 | (–, 0.50) |
|  | 30 | – | – | 1.87 | 2.09 | (–, 0.53) |
| *Krishna & Chatterjee* (2022) [29] | 10 | – | – | 3.17 | 0.62 | – |
|  | 30 | – | – | 1.95 | 2.05 | – |
| *Liu & Xu* (2021) [24] | 10 | – | – | 3.03 | 0.62 | – |
| *current study* | 10 | 1.03 | 2.05 | 3.08 | 0.62 | (0.26, 0.46) |
|  | 20 | 0.71 | 1.43 | 2.14 | 1.20 | (0.48, 0.55) |
|  | 30 | 0.43 | 1.17 | 1.60 | 1.58 | (0.53, 0.60) |

## 5.2 Results

We solved the rounded squares cylinders for $\alpha = \{2, 3\}$ up to[8] $Re = 40$, and the square cylinder (which can be formally regarded as $\alpha = \infty$) up to[9] $Re = 30$.

Figure 13 depicts the wake for $Re = 20$ and various values of $\alpha$. For the square cylinder (bottom), note how the horizontal velocity undergoes a sizeable jump discontinuity at the downstream corners. For $\alpha = 3$, the vorticity contours for several values of $Re$ are displayed in Figure 14.

The vorticity profile along the obstacle, as the latter loses smoothness, is shown in Figure 15. Note how the vorticity gradients get steeper around the rounded corners with increasing $Re$. In the square cylinder case, the RBF-PU profile breaks down in the presence of the discontinuities. (Results are similar for the pressure profile on the surface.) This reflects Gibbs-like oscillations in the underlying interpolants. The overshoot is not dampened by local refinement [18].

Finally, in Table 2 we compile the wake features of our RBF-PU simulations, and compare with the literature for the square cylinder case. (We are not aware of comparable published data for rounded cylinders.)

---
[8] We assume that the instability occurs after $Re = 40$.
[9] This upper $Re$ value is taken from the literature.

# 6 Conclusions

Even though it is a fairly general approach, the RBF-PU delivers very satisfactorily on the benchmark problem of viscous flow past a circular cylinder. In order to further improve, recipes from the literature tailored to this problem could be incorporated into the RBF-PU formulation. Particularly the pressure could be treated in a more sophisticated way, such as using staggered grids [6] or penalty terms for the lacking BC [20].

As we transition from the circular cylinder towards less blunt obstacles, the quality of the RBF-PU interpolants gradually deteriorate. Eventually, they break down at the square cylinder, owing to the boundary singularities in the solution at the corners. The extreme gradients around them induce spurious Gibbs-like oscillations in the RBF-PU near field, degrading the accuracy and convergence rates achievable with smooth solutions.While the Gibbs phenomenon has been successfully overcome in particular situations [2, 18, 19], a general recipe seems to remain an open problem—especially with nonlinear problems. (A rare example of the latter is given in [3, Example III].)

Beyond the flow past a circular cylinder, the RBF-PU methodology introduced in this paper offers—thanks to its accuracy, geometric flexibility, and efficient handling of nonlinearities—a viable alternative to solving elliptic systems of exterior BVPs with smooth solution by domain compression. The main advantages are: the ease of formulation and coding; the absence of unphysical boundary conditions and large computational domains; and the significantly smaller discretisation support required. In summary, treatment of spurious oscillations induced by boundary singularities is identified as the current major limitation.

**Funding.** FB was funded by grant 2018-T1/TIC-10914 of Madrid Regional Goverment in Spain. ASV and FB gratefully acknowledge mobility grants awarded by the Erasmus+ programme of the European Union.

**Data availability.** Data sharing not applicable to this article as no datasets were generated or analysed during the current study.

# Declarations

**Conflict of interest.** The authors declare that they have no conflict of interest.